\magnification 1200
\advance\hoffset by 1,5truecm
\advance\hsize by 2,3 truecm
\def\makefootline{\baselineskip=52pt\line{\the\footline}}
\vsize= 23 true cm
\hsize= 14 true cm
\overfullrule=0mm

\def\vmid#1{\mid\!#1\!\mid}

\headline={\hfill\tenrm\folio\hfil}
\footline={\hfill}\pageno=1

\newcount\coefftaille \newdimen\taille
\newdimen\htstrut \newdimen\wdstrut
\newdimen\ts \newdimen\tss
\def\vmid#1{\mid\!#1\!\mid}

\def\fspeciale{\textfont0=\tenrmp%
\scriptfont0=\sevenrmp%
\scriptscriptfont0=\fivermp%
\textfont1=\tenip%
\scriptfont1=\sevenip%
\scriptscriptfont1=\fiveip%
\textfont2=\tensyp%
\scriptfont2=\sevensyp%
\scriptscriptfont2=\fivesyp%
\textfont3=\tenexp%
\scriptfont3=\tenexp%
\scriptscriptfont3=\tenexp%
\textfont\itfam=\tenitp%
\textfont\bffam=\tenbfp%
\textfont\slfam=\tenbfp%
\def\it{\fam\itfam\tenitp}%
\def\bf{\fam\bffam\tenbfp}%
\def\rm{\fam0\tenrmp}%
\def\sl{\fam\slfam\tenslp}%
\normalbaselineskip=12pt%
\multiply \normalbaselineskip by \coefftaille%
\divide \normalbaselineskip by 1000%
\normalbaselines%
\abovedisplayskip=10pt plus 2pt minus 7pt%
\multiply \abovedisplayskip by \coefftaille%
\divide \abovedisplayskip by 1000%
\belowdisplayskip=7pt plus 3pt minus 4pt%
\multiply \belowdisplayskip by \coefftaille%
\divide \belowdisplayskip by 1000%
\setbox\strutbox=\hbox{\vrule height\htstrut depth\wdstrut width 0pt}%
\rm}

\def\vmid#1{\mid\!#1\!\mid}
\def\Card{\mathop{\rm Card}\nolimits}
\def\Ker{\mathop{\rm Ker}\nolimits}
\def\fle{\rightarrow}
\def\flep{\longrightarrow}

\null\vskip-1cm

\font\sc=cmcsc10

\newdimen\emm 
\def\pmb#1{\emm=0.03em\leavevmode\setbox0=\hbox{#1}
\kern0.901\emm\raise0.434\emm\copy0\kern-\wd0
\kern-0.678\emm\raise0.975\emm\copy0\kern-\wd0
\kern-0.846\emm\raise0.782\emm\copy0\kern-\wd0
\kern-0.377\emm\raise-0.000\emm\copy0\kern-\wd0
\kern0.377\emm\raise-0.782\emm\copy0\kern-\wd0
\kern0.846\emm\raise-0.975\emm\copy0\kern-\wd0
\kern0.678\emm\raise-0.434\emm\copy0\kern-\wd0
\kern\wd0\kern-0.901\emm}
\def\eps{\varepsilon}

\font\tendb=msbm10
\font\sevendb=msbm7

\newfam\dbfam
\textfont\dbfam=\tendb\scriptfont\dbfam=\sevendb\scriptscriptfont\dbfam=\sevendb
\def\db{\fam\dbfam\tendb}

\def\C{{\db C }}

\def\N{{\db N }}

\def\Q{{\db Q }}
\def\R{{\db R }}

\def\Z{{\db Z }}

\font\scgrand=cmcsc10 at 14pt

\def\Arc{\mathop{\rm Arc}\nolimits}

\def\cos{\mathop{\rm cos}\nolimits}

\def\picture #1 by #2 (#3){
  \vbox to #2{\hrule width #1 height 0pt depth 0pt
  \vfill\special{picture #3}}}
\def\scaledpicture #1 by #2 (#3 scaled #4){{
  \dimen0=#1 \dimen1=#2
  \divide\dimen0 by 1000 \multiply\dimen0 by #4
  \divide\dimen1 by 1000 \multiply\dimen1 by #4
  \picture \dimen0 by \dimen1 (#3 scaled #4)}}

\newdimen\margeg \margeg=0pt
\def\bb#1&#2&#3&#4&#5&{\par{\parindent=0pt
    \advance\margeg by 1.1truecm\leftskip=\margeg
    {\everypar{\leftskip=\margeg}\smallbreak\noindent
    \hbox to 0pt{\hss\bf [#1]~~}{\bf #2 - }#3~; {\it #4.}\par\medskip
    #5 }
\medskip}}

\newdimen\margeg \margeg=0pt
\def\bbaa#1&#2&#3&#4&#5&{\par{\parindent=0pt
    \advance\margeg by 1.1truecm\leftskip=\margeg
    {\everypar{\leftskip=\margeg}\smallbreak\noindent
    \hbox to 0pt{\hss [#1]~~}{\pmb{\sc #2} - }#3~; {\it #4.}\par\medskip
    #5 }
\medskip}}

\newdimen\margeg \margeg=0pt
\def\bba#1&#2&#3&#4&#5&{\par{\parindent=0pt
    \advance\margeg by 1.1truecm\leftskip=\margeg
    {\everypar{\leftskip=\margeg}\smallbreak\noindent
    \hbox to 0pt{\hss [#1]~~}{{\sc #2} - }#3~; {\it #4.}\par\medskip
    #5 }
\medskip}}

\def\messages#1{\immediate\write16{#1}}

\def\findem{\vrule height0pt width4pt depth4pt}

\long\def\demA#1{{\parindent=0pt\messages{debut de preuve}\smallbreak
     \advance\margeg by 2truecm \leftskip=\margeg  plus 0pt
     {\everypar{\leftskip =\margeg  plus 0pt}
              \everydisplay{\displaywidth=\hsize
              \advance\displaywidth  by -1truecm
              \displayindent= 1truecm}
     {\bf Proof } -- \enspace #1
      \hfill\findem}\bigbreak}\messages{fin de preuve}}

\def\resp{\mathop{\rm resp}\nolimits}
\def\resp.{\mathop{\rm resp.}\nolimits}


\parindent=0cm
\bigskip\bigskip
\centerline{\scgrand  Analytic Cliffordian Functions}
\bigskip
\centerline{by}
\medskip
\centerline{Guy {\sc Laville} and Eric {\sc Lehman}}
\centerline{Universit\'e de Caen - CNRS
UMR 6139}
\centerline{Laboratoire Nicolas ORESME}
\centerline{14032 Caen France}
\centerline{e-mail : guy.laville@math.unicaen.fr ; lehman@math.unicaen.fr}

\vskip1cm
\hskip-0,2cm{\bf Abstract}
\bigskip
In classical function theory, a function is holomorphic if and only if it is complex
analytic. For higher dimensional spaces it is natural to work in the context of
Clifford algebras. The structures of these algebras depend on the parity of the
dimension \ $n$ \ of the underlying vector space. The theory of holomorphic
Cliffordian functions reflects this dependence.
\ In the  case of odd $n$ the space of functions is defined by an operator (the
Cauchy-Riemann equation) but not in the case of even $n$. For all dimensions the powers
of identity
\
$(z^n, x^n)$ \  are the foundation of function theory.

\medskip
AMS \ classificatioon : 30 G 35, 15 A 66.

\vskip1cm
\hskip-0,5cm{\bf I. Introduction}
\bigskip
A complex analytic function\ \ $f(z)$\ \ may be defined as being
locally the sum of a convergent power series\ \ $f(z)
=\displaystyle\sum_{N=1}^\infty a_Nz^{N-1}$\ \ or as being holomorphic, that is
such that \  $\Bigl(\displaystyle{\partial\over\partial x}+i\
{\partial\over
\partial y}\Bigr)f(x+iy)= 0$. A real analytic function\ $u(x)$\ may be defined
as being locally the sum of a convergent power series or by\ $u(x) = f(x+i0)$\
where\ $f$\ is a complex holomorphic function such that\ $f(\overline z) =
\overline{f(z)}$. The main difference between the complex and the real
case is the existence or non existence of a differential relation
characterizing holomorphy.
\bigskip
We extend the definitions of analyticity and holomorphy to functions defined on
Clifford algebras\ $\R_{0,n}$\ distinguishing between the case\ of odd $n$,
$n = 2m+1$, and the case\ of even $n$, $n = 2m$. We show that the equivalence
between analyticity and holomorphy still holds. The cases\ of odd $n$
and even \  $n$ \  interrelate in a way that reflects the difference
between the structures of the algebras\ \ $\R_{0,2m}$\ and\
$\R_{0,2m+1}$. In particular the center of\ $\R_{0,2m}$\
is\ $\R$\ although the center of\ $\R_{0,2m+1}$\ is \ $\R\oplus\R
e_{12\ldots 2m+1}$, where\ $e_{12\ldots 2m+1}$\ is a pseudoscalar.

\vskip1cm
\hskip-0,7cm{\bf II. Notations}
\bigskip 
Let\ $V_n$\ be an anti-Euclidean vector space of dimension\ $n$. For any
orthonormal basis\ \ $e_1,\ldots ,e_n$\ \ of\ \ $V_n$\ \ we have for
all distinct\ \ $i$\ \ and\ \ $j$\ \ in\ \ $\{ 1,\ldots ,n\}$ 

$$e_i^2 =-1\qquad\hbox{and}\qquad e_ie_j = -e_je_i.$$

If\ \ $I\subset\{ 1,\ldots ,n\}$\ \ and\ \ $I = \{i_1,\ldots ,i_k\}$\ \
with\ \ $i_1<\ldots <i_k$\ \ we set\ \ $e_I = e_{i_1}e_{i_2}\ldots
e_{i_k}$. For\ \ $I=\emptyset$, we set\ \ $e_\emptyset = e_0 = 1$. Then\
\ $(e_I)_{I\subset\{1,\ldots ,n\}}$\ is a basis of the Clifford
algebra\ \ $\R_{0,n}$\ \ seen

\bigskip
 as a real vector space. If\ \ $A =
\displaystyle\sum_{I\subset\{1,\ldots ,n\}} A_Ie_I$, with\ \ $A_I\in\R$, is an
element of\ \ $\R_{0,n}$
 \ we call\ \ $A_0=A_\emptyset$\ \ the scalar
part of\ \ $A$\ \ and denote it by\ \ $A_0 = S(A)$. Following the \ $\R - \C$
\ case and also Leutwiler and Eriksson-Bique\ [EL2], we introduce the
decomposition

$$\R_{0,n} =\R_{0,n-1}\oplus e_n\R_{0,n-1}$$
\bigskip
(For convenience in our computations we have chosen\ \
$e_n\R_{0,n-1}$\ \ instead of\ \ $\R_{0,n-1}e_n$). This decomposition
means that given a vector\ \ $e_n$\ \ there are two maps from\ $\R_{0,n}$\ to\
\ $\R_{0,n-1}$, denoted\ \ ${\cal R}$\ and\ ${\cal J}$, such that for
any\ \ $A$\ \ in\ \ $\R_{0,n}$\ \ we have

$$A = {\cal R}A +e_n {\cal J}A.$$

We have chosen  notation \ \ ${\cal R}$\ \ and\ \ ${\cal J}$\ and the following
notation for conjugation to stress the fact that for \ $n=1$, \ we have \ 
$\R_{0,1} = \C$, $\R_{0,0} = \R$ \ which yield the usual relations between\ \ $\C$
\ and
\
$\R$. If\
\
$A\in\R_{0,n}$, we call the conjugate of \ $A$ \ and denote by \ $\overline A$\ \
the element of\ \
$\R_{0,n}$\ \ defined by 

$$\overline A = {\cal R} A-e_n{\cal J} A$$ 

If\ $z$\ is a paravector, that is an element of\ $\R\oplus V_n$, we have

$$z = z_0+z_1e_1 +\ldots +z_ne_n\qquad \hbox{with}\qquad z_0\in\R,
z_1\in \R,\ldots ,z_n\in\R.$$

We denote by \ $|z|$ \ the positive real number such that \
$|z|^2=z_0^2+z_1^2+\ldots +z_n^2$ \ and by \ $x = {\cal R}z =
z_0+z_1e_1+\ldots +z_{n-1}e_{n-1}$ \ the paravector in\ \ $\R\oplus
V_{n-1}$\ \ such that 
$$z = x+z_ne_n\qquad \hbox{and}\qquad \overline z = x-z_ne_n.$$

We define \ $z_\ast$ \ by : $z_\ast = z_0-z_1e_1-\ldots -z_ne_n$. Then\
\ $|z|^2 = zz_\ast = z_\ast z.$ 
\medskip
We introduce the  differential operators 

$$D = {\partial\over \partial z_0} +e_1 {\partial\over \partial z_1}
+\ldots +e_n{\partial \over\partial z_n},\ \ D_\ast = {\partial\over
\partial z_0} - e_1 {\partial\over \partial z_1} -\ldots -e_n{\partial
\over\partial z_n}$$ and $$\Delta = DD_\ast = D_\ast D.$$

Note that \ $\Delta$ \ is the usal Laplacian.
\medskip
If\ \ $\alpha = (\alpha_0,\alpha_1,\ldots ,\alpha_n)\in (\N\cup \{
0\})^{n+1}$\ \ is a multi-index we denote its length by\ \ $|\alpha | =
\alpha_0+\alpha_1+\ldots +\alpha_n$. The elementary multi-index\ \
$\eps_k$\ \ is defined by\ \ $\eps_k =
(\delta_{k0},\delta_{k1},\ldots ,\delta_{kn})$\ \ where\ \
$\delta_{ij}$\ \ is the Kronecker symbol equal  to\ \ 1\ \ if\ \
$i=j$\  \ and to\ 0\ if\ \ $i\not= j$.

The order on the set of multi-indexes is the lexicographical order.

\vskip1cm
\hskip-1cm{\bf III. Analytic Cliffordian polynomials}
\bigskip 
Laville and Ramadanoff \ [LR1] \ have defined holomorphic Cliffordian polynomials
for odd  \ $n$. The same definitions can also be used for even \ $n$. We will see that
these homogeneous polynomials are the building blocks of analytic Cliffordian functions
in both cases. Therefore we call them analytic Cliffordian polynomials. Polynomials
with the same structure were introduced by Heinz Leutwiler in [Le1] and also in [EL1].
Monogenic polynomials are particular cases of analytic Cliffordian polynomials [BDS],
[DSS].

\bigskip
{\bf III.1. Three classes of analytic Cliffordian polynomials}
\medskip
{\bf Definition 1.a.-} Let \ $a$ \ be a paravector and\ \ $N\in\N$.
We call elementary analytic monomial and denote by \ $M_N^a(z)$ \ the
homogeneous monomial function of degree \ $N-1$ \ defined by 

$$M_N^a(z) = (az)^{N-1} a = a(za)^{N-1}.$$

{\bf Remark 1.-} \  We may also define\ \ $M_N^a(z)$\ \ for all integers
by\ \ $M_0^a(z) = z^{-1}$\ \ and for\ \ $N < 0$, $M_N^a(z) =
M_{1-N}^{z^{-1}}(a^{-1})$. It is often convenient to write\ \
$a\kern-0,18cm(\,z\,a\kern-0,18cm)^{N-1}$\ instead of\ $M_N^a(z)$. We have
then for\ $N\in \Z$\ and\ $\sqrt a$\ a paravector such that\ $(\sqrt a)^2 = a$

$$a\kern-0,18cm(\,z\,a\kern-0,18cm)^N  = \sqrt a(\sqrt az\sqrt a)^N\sqrt
a.$$
These polynomials are close to similar ones in [Le2].
\bigskip
{\bf Remark 2.-} \ $M_N^a(z) = |a|^N M_N^{a/|a|}(z) = |z|^{N-1}
M_N^a(z/|z|).$ 

\bigskip
{\bf Proposition 1.-} \  $M_N^a(z)\in\R\oplus V_n$\ \ and\ \
$|M_N^a(z)| = |a|^N|z|^{N-1}$.
\bigskip
{\bf Proof.-}
Computing \ $aza$ \ explicitly, we get \ $aza\in\R\oplus V_n$ \ and \
$|aza| = |a|^2|z|$. Note that 
$$M_N^a(z) = aM_{N-1}^z (a)a.\leqno (*)$$

\bigskip
{\bf Proposition 2.-}\ \ $M_N^a(z)= D_\ast \displaystyle{1\over N} S((az)^N).$
\bigskip
{\bf Proof.-}
Let us define\ \ $\theta$\ \ by\ \ $|a||z|\cos \theta = S(az),\ 0\leq\theta
\leq\pi$. Then \ $2S(az) = az+z_*a_*$ \ and  \ $aza = (az+z_*a_*)a-z_*a_*a =
2S(az)a - \vmid{a}^2 z_* = 2 \vmid{a} \ \vmid{z} \cos \theta a - \vmid{a}^2
\vmid{z}^2 z^{-1}$. A simple recursion using (*) yields 
$$M_N^a(z) = |a|^{N-1} |z|^{N-1} {\sin N\theta\over \sin \theta} a -
|a|^N|z|^N {\sin (N-1)\theta \over \sin\theta} z^{-1}.$$ 
Since \ $2S((az)^N) = M_N^a(z) z + z_* \ M_N^a(z)_*$, \ we get
$$S\bigl((az)^N\bigr) = |a|^N |z|^N \cos (N\theta ).$$

From the definitions of \ $|z|$ \ and \ $\theta$, \ we get\ \ $D_\ast|z| =
|z|z^{-1}$ and  $D_\ast\theta = \displaystyle{\cos \theta\over
\sin\theta} z^{-1} - |a|^{-1}|z|^{-1}\displaystyle{1\over \sin\theta} a$.
Finally

$$D_\ast S(az)^N = N|a|^N|z|^N z^{-1} \cos (N\theta )-N|a|^N|z|^N \sin N\theta
D_\ast \theta = NM_N^a(z).$$

{\bf Corollary.-}\ \ $z^N = D_\ast \displaystyle{1\over N+1} S(z^{N+1}).$
\medskip
{\bf Proof.-}
Chose\ \ $a = e_0$\ \ and replace\ \ $N$\ \ by\ \ $N+1$\ \ in proposition 2.

\bigskip
{\bf Remark.-} \ Let \ $T_N(x)$ \ and \ $U_N(x)$, $x\in\R$, be  the classical
Tchebycheff polynomials of the first and second kind. Recall that
$$\eqalign{
&T_N(x) = \cos~\!(N~\Arc \cos x)\cr
&U_N(x) = {\sin~\! \bigl( (N+1) \Arc \cos x\bigr)\over \sin~\!(\Arc \cos x)}.\cr}$$

Thus, when \ $\vmid{a} = 1$, \ $\vmid{z} = 1$, we get  
$$S\bigl( (az)^N\bigr) = T_N\bigl( S(az)\bigr) \ \hbox{and} \ M_N^a (z) = U_{N-1}
\bigl( S(az)\bigr) a - U_{N-2} \bigl( S(az)\bigr) z^{-1}.$$

\bigskip
{\bf Definition 1.b.-} Let\ \ $a_1,\ldots ,a_k$\ \ be paravectors and let \ 
$N_1,\ldots ,N_k$ \ be  integers belonging to\ \ $\N\cup\{0\}.$ We set\ \ $N =
N_1+\ldots +N_k$\ \ and denote by\ \ ${\cal P}_{N_1,\ldots ,N_k}$\ \ the set
of partitions\ $I$\ of\ $\{1,\ldots ,N\}$\ into a union of disjoint
subsets\ $I = (I_1,\ldots ,I_k)$\ \ such that\ $\Card I_1 = N_1,\ldots
,\Card I_k = N_k$. For\ $I\in{\cal P}_{N_1,\ldots , N_k}$\ and\ $\nu \in\{
1,\ldots ,N\}$, we define\ $b_\nu^I$\ by\ $b_\nu^I = a_j$\ where\ $j$\ is the
element of\ $\{1,\ldots ,k\}$\ such that\ $\nu\in I_j$. We define the\
$a_1,\ldots ,a_k$\ symmetrical analytic homogeneous polynomial of degree\ $N-1$\
in\ $z$\ and\ $N_j$\ in\ $a_j$\ by 

$$S_{N_1,\ldots ,N_k}^{a_1,\ldots ,a_k} (z) = \sum_{I\in{\cal P}_{N_1,\ldots
,N_k}} \Bigl(\prod_{\nu=1}^{N-1}  (b_\nu^I z)\Bigr) b_N^I.$$

\bigskip
{\bf Proposition 3.-}\ \ $S_{N_1,\ldots ,N_k}^{a_1,\ldots ,a_k}$\ \ is a real
linear combination of elementary analytic Cliffordian monomials of degree\
$N_1+\ldots +N_k-1.$

\medskip
{\bf Proof.-} For any real \ $\lambda$, we have
$$M_N^{a+\lambda b} (z) = \sum_{p+q=N} \lambda^q S_{p,q}^{a,b} (z).$$
Choose\ $N+1$\ different values for\ $\lambda$; one gets a Van der Monde matrix
which is invertible. This shows the result for\ $k=2$. We can iterate the same
argument noting that for any real\ $\lambda$

$$S_{N_1,\ldots ,N_{k-1},N_k}^{a_1,\ldots ,a_{k-1},a_k+\lambda a_{k+1}}(z) =
\sum_{p+q=N_k} \lambda^q S_{N_1,\ldots ,N_{k-1},p,q}^{a_1,\ldots
,a_{k-1},a_k,a_{k+1}}(z).$$

\bigskip\medskip
{\bf Definition 2.-} For each multi-index\ $\alpha\in (\N\cup \{ 0\})^{n+1}$\ we
define  \ $Q_\alpha$ \ as the homogeneous polynomial in\ $z_0,\ldots ,z_n$\ of degree\
$|\alpha |-1$,  by 
$$Q_\alpha (z) = \partial_\alpha z^{2|\alpha |-1}$$
where\ $\partial_\alpha$\ is the differential operator of order\ $|\alpha | :
\partial_\alpha = \displaystyle{\partial^{|\alpha |}\over \partial
z_0^{\alpha_0}\partial z_1^{\alpha_1}\ldots \partial z_n^{\alpha_n}}.$

\bigskip
{\bf Definition 3.-} For each multi-index\ $\alpha\not= (0,0,\ldots ,0)$\ we
define the analytic Cliffordian polynomial\ $P_\alpha$\ by 
$$P_\alpha (z) = \sum_{\sigma\in S_\alpha}\Bigl(\prod_{\nu =1}^{|\alpha
|-1}(e_{\sigma (\nu )}z)\Bigr) e_{\sigma (|\alpha |)}$$

where\ $S_\alpha$\ is the set of maps\ $\sigma$\ from\ $\{ 1,\ldots ,|\alpha |\}$\
to\ $\{ 0,1,\ldots ,n\}$\ such that\ $\Card (\sigma^{-1}(\{k\})) = \alpha_k$\ for
all\ $k$\ in\ $\{ 0,1,\ldots ,n\}$.

\bigskip\bigskip
{\bf III.2. Relations among the\ $M_N^a(z)$, the\ $Q_\alpha (z)$\ and the\
$P_\alpha (z)$}
\bigskip
For any \  $N$ \ in \ $\N$, the real linear space generated by the elementary analytic
monomials\ $M_N^a$, the real linear space generated by the\ $Q_\alpha$\ with\
$|\alpha | = N$\ and the real linear space generated by the\ $P_\alpha$\  with\
$|\alpha | = N$\ are identical. 

\bigskip
{\bf Proposition 4.-}\ \ $Q_\alpha (z) = \partial_\alpha M_{2|\alpha |}^1 (z).$
\medskip
{\bf Corollary.-} For any multi-index\ $\alpha$, $Q_\alpha (z) \in \R\oplus V_n$ \
and there exists a scalar polynomial\
$q_\alpha (z)$\ homogeneous of degree\ $|\alpha |$\ such that\ $Q_\alpha (z) =
D_\ast q_\alpha (z)$.

{\bf Proof.-}
$\partial_\alpha$\ and\ $D_\ast$\ commute and\ $M_{2|\alpha |}^1|z| = z^{2|\alpha
|-1}$. Use proposition 1 and corollary of proposition 2.

\bigskip
{\bf Proposition 5.-}  $$Q_\alpha (z) = k_\alpha P_\alpha(z) +
\displaystyle\sum_{\alpha'>\alpha\atop |\alpha'| = |\alpha |}
\lambda_{\alpha\alpha'} P_{\alpha'}(z)$$ where $$k_\alpha = \pmatrix{2|\alpha
|-1\cr \alpha_0\cr}\ \alpha_0!\alpha_1!\ldots \alpha_n!$$ and
$$\lambda_{\alpha\alpha'}\in\Z.$$

\bigskip
Proof.- \  Let \ $\beta = (\alpha_1,\ldots,\alpha_n) \in (\N -
\{ 0\})^n$ \ be the multi-index such that \ $\alpha =
(\alpha_0,\beta)$. From the definition of \ $Q_\alpha$ \
follows~:
$$Q_\alpha (z) = \pmatrix{
2\vmid{~\!\!\alpha~\!\!} -1\cr
\alpha_0\cr} \alpha_0 ! \ {\partial^{\vmid{~\!\beta~\!}} \over
\partial z_1^{\alpha_1} \ldots \partial z_n^{\alpha_n}} \ 
z^{2\vmid{~\!\beta~\!} -1+\alpha_0}.$$

\bigskip
Consider \ $z^{2\vmid{~\! \beta~\!} - 1 + \alpha_0}$ \ as an
explicit product of \ $2\vmid{~\! \beta~\!}-1+\alpha_0$ \
factors each equal to $z$, that is~: \ $z^{2\vmid{~Ê\! \beta~Ê\!
}-1+\alpha_0} = z \cdot z \cdot z\cdot \ \cdots z$.  To apply \
$\displaystyle{\partial\over\partial x_k}$ \ is equivalent to
replace each \ $z$ \ once by \ $e_k$ \ and to add all the
products obtained. If we never derivate two successive \ $z$
\ then we get \ $\alpha_1 ! \ldots \alpha_n ! \ P_\alpha(z)$.
If we derivate two successive \ $z$, \ we get factors \ $e_k
e_h$ \ which anihilate if \ $k \not= h$ \ because \ $e_ke_h =
- e_h e_k$ \ and else \ $-1$ \ if \ $k = h$. Since \ $1
= e_0$ \ we get the terms of \ $P_{\alpha'}$ \ with \ $\alpha'
> \alpha$ \ in the lexicographical order.

\vskip1cm
{\bf Corollary 1.-} \ $P_\alpha (z) = k_\alpha^{-1} Q_\alpha (z)
+\displaystyle\sum_{\alpha'>\alpha\atop |\alpha'|=|\alpha
|}\mu_{\alpha\alpha'} Q_{\alpha'}(z)$ \ where \ $\mu_{\alpha\alpha'}\in \Q$.

\medskip
{\bf Corollary 2.-} For any multi-index \ $\alpha$, \ $P_\alpha (z) \in \R\oplus V_n$ 
\ and there exists a scalar polynomial \ $p_\alpha (z)$ \ homogeneous of degree \
$|\alpha |$\  such that \ $P_\alpha (z) = D_\ast p_\alpha (z)$.

\bigskip\medskip
{\bf Proposition 6.-}\ \ $M_N^a(z) = \displaystyle\sum_{|\alpha | = N}
a^\alpha P_\alpha (z)$, where\ $a^\alpha :=
a_0^{\alpha_0}a_1^{\alpha_1}\ldots a_n^{\alpha_n}.$

\bigskip\medskip
{\bf Proposition 7.-} \  For\ $\alpha = (\alpha_0,\alpha_1,\ldots ,\alpha_n)$, we have
\  $P_\alpha (z) = S_{\alpha_0,\alpha_1,\ldots ,\alpha_n}^{e_0,e_1,\ldots ,e_n}
(z).$
\medskip
Proof.- \ The definiton of \ $P_\alpha$ \ is the definition of
\ $S$ \ in which the \ $k$  paravectors \ $a_1 \cdots a_k$
\ are the \ $(n+1)$ \ elements of a basis of  $S\oplus V :
e_0, e_1,\ldots,e_n$.

\bigskip\medskip
{\bf Corollary.-} The real linear space generated by the\ $P_\alpha$\ with\
$|\alpha | = N$\ is independent of the basis\ $e_1,\ldots ,e_n$\ of\ $V_n$.
\medskip
{\bf Remarks.-} \  In the case of odd \ $n$, \  this is already known. 

Note that generally the polynomials \ $P_\alpha (z)$\ with \ $|\alpha | = N$\ are not \
$\R$-linearly independant. For example, if\ $n=3$ and \ $N=4$ \ we have 
$$3P_{4000} + 3P_{0400} +3P_{0040} +3P_{0004} +P_{2200} + P_{2020} +
P_{2002} +P_{0220}+ P_{0202} +P_{0022} = 0.$$

\bigskip\medskip
{\bf  III.3. Some properties of the polynomials \  $P_\alpha$ }
\bigskip
{\bf Proposition 8.-} For\ $\alpha = (\alpha_0,\ldots ,\alpha_n)$, the polynomial \
$p_\alpha (z)$\ is even in\ $z_k$\ if \ $\alpha_k$ \ is even and   odd in\ $z_k$ \ if \
$\alpha_k$\ is odd.

{\bf Proof.-}
We know that\ $P_\alpha (z)\in\R\oplus V_n$. Then we can write 
$$P_\alpha (z) = A_0(z) + A_1(z)e_1+\ldots +A_n(z)e_n.$$

Suppose\ $\alpha_k$\ even.  $P_\alpha (z)$\ is a sum of terms like

$$e_{i_1}ze_{i_2} z\ldots ze_{i_{|\alpha |}}$$
in which\ $e_k$\ occurs\ $\alpha_k$\ times, that is an even number of times.
Write\ $z = z_0+z_1e_1+\ldots +z_ne_n$\ and develop all the terms. The terms
 which contribute to\ $A_k(z)e_k$\ must contain an odd number of times
the vector\ $e_k$\ and then\ $z_ke_k$\ has to appear an odd number of
times. $A_k(z)$\ is then a sum of terms which are all odd in\ $z_k$\ and\
$A_k(z)$\ is odd in\ $z_k$. Since\ $P_\alpha (z) = D_\ast p_\alpha (z)$, we
have\ $A_k(z) = - \displaystyle{\partial\over \partial z_k}\ p_\alpha (z)$.
Since\ $p_\alpha (z)$\ is homogeneous, we can conclude that\ $p_\alpha (z)$\
is even in\ $z_k$.
\medskip
If\ $\alpha_k$\ is odd the proof is the same.

\medskip
{\bf Examples :} \ $p_{(2,0,1,0)}(z) = (3z_0^2-z_1^2-z_2^2-z_3^2)z_2$ \ is even
in\ $z_0,z_1$\ and\ $z_3$\ and odd in\ $z_2$.\  $p_{(1,1,1,0)}(z)
=8z_0z_1z_2$\ is odd in\ $z_0,z_1$ \ and\ $z_2$\ and even in\ $z_3$.

\medskip
{\bf Corollary.-} If\ $\alpha_n$\ is even, then\ $P_\alpha(\overline z) =
\overline{P_\alpha(z)}$ ; if\ $\alpha_n$\ is odd then\ $P_\alpha(\overline z)
= -\overline{P_\alpha (z)}$.
\medskip
{\bf Proof.-}
Let us write
$$P_\alpha (z) = A_0(z) + A_1(z)e_1+\ldots + A_{n-1}(z)e_{n-1} + A_n(z)e_n.$$

If \ $\alpha_n$ \ is even, then\ $p_\alpha (z)$\ is even in\ $z_n$\ and the
polynomials \
$A_k =
\pm \displaystyle{\partial\over\partial z_k} p_\alpha (z)$\ for\ $k\not= n$,
are all even in\ $z_n$, but\ $A_n$\ is odd in\ $z_n$. Thus

$$P_\alpha (\overline z) = A_0(z) + A_1(z)e_1 +\ldots +A_{n-1}(z)e_{n-1} -
A_n(z)e_n = \overline{P_\alpha (z)}.$$

If\ $\alpha_n$\ is odd,\ $A_0,A_1\ldots A_{n-1}$\ are odd in\ $z_n$\ and\
$A_n$\ is even in\ $z_n$, so that

$$P_\alpha (\overline z) = -A_0(z) - A_1(z)e_1-\ldots - A_{n-1}(z)e_{n-1}
+A_n(z)e_n = -\overline{P_\alpha (z)}.$$

\bigskip
{\bf Remark.-}\ $\displaystyle{\partial\over \partial z_i} P_\alpha (z)
= e_i^2\Bigl\{ 2|\alpha |P_{\alpha -\eps_i} (z) - (\alpha_i+1)\sum_{k=0}^n
P_{\alpha +\eps_i-2\eps_k}(z)\Bigr\}$.

\vskip1cm
\hskip-0,9cm{\bf IV. Analytic Cliffordian functions and holomorphic
Cliffordian functions}

\bigskip
{\bf Definition 1.-} Let\ $\Omega$\ be a domain of\ $\R\oplus V_n$\ and\ $f :
\Omega\fle \R_{0,n}$. We say that\ $f$\ is a left analytic Cliffordian
function if  any\ $\omega$ \ in \ $\R\oplus V_n$ \ has a neighbourhood\
$\Omega_\omega$ \ in\ $\Omega$\ such that for any\ $z$\ in\
$\Omega_\omega$,\ $f(z)$\ is the sum of a convergent series
$$f(z) = \sum_{N=1}^\infty\sum_{a\in A_N} M_N^a(z-\omega )C_a$$

where for each\ $N$\ in\ $\N$,\ $A_N$\ is a finite subset of\ $\R\oplus
V_n$, for each a in\ $A_N$, $C_a\in\R_{0,n}$\ and\
$\displaystyle\sum_{N=1}^\infty\displaystyle\sum_{a\in A_N} |a |^N|z-\omega
|^{N-1}|C_a|$\ is convergent in \ $\Omega_\omega$.

\bigskip
{\bf Remark 1.-} \  The relation\ $M_N^a(z-\omega ) = \displaystyle\sum_{p+q=N}
(-1)^q S_{p,q}^{a,a\omega a}(z)$\ and the proposition 3 prove the consistency
of the definition with respect to translations. Consequentely we will restrict
ourselves to the case\ $\omega = 0$.
\medskip
{\bf Remark 2.-} \  The above definition is obviously intrinsic, but we get an
equivalent definition if we replace the \ monomials \  $M_N^a$\ by the polynomials \
$P_\alpha$; the function\ $f : \R\oplus V_n \fle\R_{0,n}$\ is left analytic
Cliffordian in a neighbourhood\ $\Omega$\ of\ $0$\ if for every\ $z$\ in\ $\Omega$,
$f(z)$\ is the sum of a convergent series 

$$f(z) = \sum_{N=1}^\infty\sum_{|\alpha |= N} P_\alpha (z) c_\alpha$$

where\ $\alpha$\ are multi-indexes belonging to\ $(\{ 0\}\cup\N)^{1+n}$, and
for each\ $\alpha$\ we have\ $c_\alpha\in\R_{0,n}$\ and\
$\displaystyle\sum_{N=1}^\infty\sum_{|\alpha |=N} |P_\alpha (z)|\ |c_\alpha|$\
is convergent.

\bigskip
{\bf Definition 2.-} Let\ $\Omega$\ be a domain in\ $\R\oplus V_n$.  A
fonction\ $u : \R\oplus V_n\fle \R_{0,n}$\ is called a left holomorphic
Cliffordian function 

{\parindent =1cm
\item{\bf (i)} for odd \ $n$,  if\ $D\Delta^mu = 0$ \ where \ $m = (n-1)/2$
\item{\bf (ii)} for even\ $n$, \  if for any\ $\omega\in\Omega$\  a neighbourhood\
$\Lambda_\omega$\ in\ $\R\oplus V_{n+1}$\ and a left holomorphic Cliffordian function\
$f$\ defined on \ $\Lambda_\omega$ \ exist such that 

\hskip2cm - for all\ $z$\ in \ $\Lambda_\omega$ : \ $\bar z$ is in $\Lambda_\omega$ and
$f(\overline z) =
\overline{f(z)}$

\hskip2cm - for all\ $x$\ in\ $\Lambda_\omega \cap (\R\oplus V_n)$,\ $u(x) =
f(x)$.
\par}
\medskip
\input epsf.def
$$\epsfbox{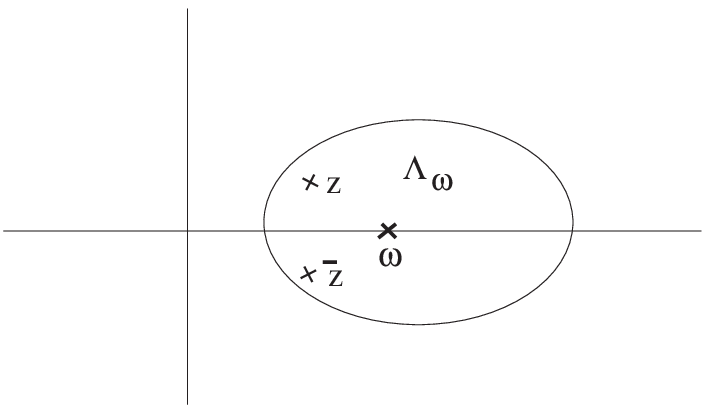}$$ 
\vskip-1cm
\hskip3,5cm$\R\oplus V_{n+1}$ 

\vskip-4,cm
\hskip5,5cm $\R e_n$
\vskip2cm
\hskip11cm $\R\oplus V_n$
\bigskip\bigskip

\vskip1cm
{\bf Theorem.-} Let \ $\Omega$\ be a domain of\ $\R\oplus V_n$. A function\ $f
: \Omega\fle \R_{0,n}$\ is left analytic Cliffordian if and only if it is left
holomorphic Cliffordian.
\medskip
For odd \   $n$, the theorem has already been proven in\ [LR1]. Let\ $n$\ be even, $n
= 2m$.  
Let \ $u$ \ be a left analytic Cliffordian function on a neighbourhood\
$\Omega$\ of 0 and let\ $S_n(r)$\ be a sphere of center 0 and radius\ $r >
0$\ included in\ $\Omega$. For\ $|x| < r$\ we have

$$u(x) = \sum_{N=1}^\infty \sum_{a\in A_N} M_N^a(x)\ C_a$$

where\ $\displaystyle\sum_{N=1}^\infty \sum_{a\in A_N} |a|^N |x|^{N-1} |C_a|$\
is convergent. Chose\ $\Lambda_0 = S_{n+1}(r)$\ the interior of the sphere of center
0 and radius\ $r$\ in\ $\R\oplus V_{n+1}$\ and let\ $f : S_{n+1}(r) \fle \R_{0,n+1}$\
be defined by
$$f(z) = \sum_{N=1}^\infty \sum_{a\in A_N} M_N^a (z) \ C_a.$$  

Then\ $f$\ is left holomorphic Cliffordian on\ $S_{n+1}(r)$\ and\ $f(x) =
u(x)$. Since\ $a\in \R\oplus V_n$, one gets\ $M_N^a(\overline z) =
\overline{M_N^a(z)}$. And since\ $C_a\in\R\oplus V_n$, we have\ $f(\overline z
) = \overline{f(z)}$.
\medskip
Conversely, let\ $f$\ be a left holomorphic Cliffordian function defined on a
neighbourhood\ $\Lambda_0$\ of 0 in\ $\R\oplus V_{n+1}$\ such that\
$f(\overline z) =\overline{f(z)}$. We want to show that\ $u : \Lambda_0\cap
\R\oplus V_n\fle \R_{0,n}$,\ $x\longmapsto u(x) = f(x)$\ is left analytic
Cliffordian.
\medskip
Since\ $n+1$\ is odd, $f$ is analytic Cliffordian and we can write 
$$f(z) = \sum_{N=1}^\infty \sum_{|\beta | = N} P_\beta (z) c_\beta = D_\ast
\sum_{N=1}^\infty \sum_{|\beta | = N} p_\beta (z)c_\beta.    
$$

Let\ $H_N$\ be the real linear space of scalar homogeneous polynomials in\
$z_0$,$\ldots,$ $z_n$, $z_{n+1}$\ of total degree\ $N$\ and\ $(m+1)$-harmonic
generated by\ $(p_\beta )_{|\beta |=N}$\ (the dimension of\ $H_N$\ is\
$C_{N+n}^n-C_N^n$). We can extract a subset\ $B_N$\ of\ $\{\beta \in (\{
0\}\cup\N)^{n+2} / |\beta | = N\}$\ such that\ $(p_\beta )_{\beta\in B_N}$\ is
a basis of\ $H_N$. $D_\ast$\ is a linear map from\ $H_N$\ to\ $\R\otimes
V_{n+1}$\ and\ $\Ker D_\ast$ is a subspace of\ $H_N$. Let\
$(\psi_j)_{j\in{\cal J}}$\ be a basis of\ $\Ker D_\ast$ ; since\ $(p_\beta
)_{\beta\in B_N}$\ is a basis of\ $H_N$\ there is a subset\ $B_N^\circ$\ of\
$B_N$\ such that\ $((\psi_j)_{j\in {\cal J}}$, $(p_\beta)_{\beta\in
B_N^\circ})$\ is a basis of\ $H_N$. Then\ $((\psi_j\otimes e_I)_{j\in{\cal
J},I\subset\{ 1,\ldots ,n+1\}}$, $(p_\beta\otimes e_I)_{\beta\in B_N^\circ
,I\subset\{ 1,\ldots ,n+1\}}$\ is a basis of the real linear space\
$H_N\otimes\R_{0,n+1}$, and there are unique real numbers\ $\theta_{j,I}$\ and\
$d_{\beta ,I}$\ such that\ $\displaystyle\sum_{|\beta |=N} p_\beta (z)c_\beta
= \sum_{j\in{\cal J}}\sum_{I\subset\{ 1,\ldots ,n+1\}}\
\theta_{j,I} \psi_j(z)e_I +\sum_{\beta\in B_N^\circ}
\sum_{I\subset\{1,\ldots ,n+1\}} d_{\beta ,I} p_\beta (z) e_I$. Let us write\ $d_\beta
= \sum_I d_{\beta ,I}e_I$, using homogeneity we get for any analytic cliffordian
function\ $f$\ the existence and unicity of the coefficients\ $d_\beta $\ in\
$\R_{0,n+1}$\ such that 

$$f(z) = \sum_{N=1}^\infty\sum_{\beta\in B_N^\circ} P_\beta (z)d_\beta.$$

Let\ $B_N^{\circ+}$\ be the set of multi-indexes \ $\beta = (\beta_0,\ldots
,\beta_n, \beta_{n+1})$\ where\ $\beta \in B_N^\circ$\ and\ $\beta_{n+1}$\ is
even and\ $B_N^{\circ -} = B_N^\circ - B_N^{\circ +}$. If\ $\beta\in
B_N^{\circ +}$,  the corollary of proposition 8 implies \ $P_\beta
(\overline z) = \overline{P_\beta (z)}$\ and if\ $\beta \in B_N^{\circ -}$\ we
have\  $P_\beta (\overline z) = -\overline{P_\beta (z)}$. The relation\
$f(\overline z) = \overline{f(z)}$\ becomes then  
$$\sum_{N=1}^\infty \Bigl\{\sum_{\beta\in B_N^{\circ +}} \overline{P_\beta
(z)}~d_\beta - \sum_{\beta \in B_N^{\circ -}} \overline{P_\beta (z)}~
d_\beta\Bigr\} = \sum_{N=1}^\infty \Bigl\{\sum_{\beta\in B_N^{\circ
+}}\overline{P_\beta (z)}~\overline d_\beta +\sum_{\beta \in B_N^{\circ
-}}\overline{P_\beta (g)}~\overline d_\beta\Bigr\}$$

By homogeneity we get 
$$\sum_{\beta \in B_N^{\circ +}} \overline{P_\beta (z)}(\overline d_\beta
-d_\beta ) + \sum_{\beta \in B_N^{\circ -}}\overline{P_\beta (z)} (\overline
d_\beta + d_\beta ) = 0$$

and by conjugation 
$$\sum_{\beta \in B_N^{\circ +}} P_\beta (z) (d_\beta -\overline d_\beta )
+\sum_{\beta \in B_N^{\circ -}} P_\beta (z) (d_\beta +\overline d_\beta ) = 0.$$

By unicity of the coefficients of the\ $P_\beta$\ for\ $\beta \in B_N^\circ$,
we get\ $d_\beta = \overline d_\beta$\ if\ $\beta\in B_N^{\circ +}$\ and\
$\overline d_\beta = -d_\beta$\ if\ $\beta\in B_N^{\circ -}$. Let us write\
$a_\beta = d_\beta = d_{\overline \beta}$\ if\ $\beta\in B_N^{\circ +}$\ and\ $e_{n+1}
b_\beta = d_\beta = -\overline d_\beta$\ if\ $\beta\in B_N^{\circ -}$, we
 get
$$f(z) = \sum_{N=1}^\infty \left\{   \sum_{\beta\in B_N^{\circ +}}
P_\beta (z) a_\beta + \sum_{\beta \in B_N^{\circ -}} P_\beta (z) e_{n+1}
b_\beta
\right\},$$

where\ $a_\beta \in\R_{0,n}$\ and\ $b_\beta\in \R_{0,n}$. The two
following lemmas will then prove the  theorem.

\bigskip
{\bf Lemma 1.-} Let\ $n = 2m$\ and\ $\beta = (\beta_0,\ldots ,\beta_n,
\beta_{n+1})$. If\ $\beta_{n+1}$\ is even then the restriction of\ $P_\beta$\
to\ $\R\oplus V_n$\ is analytic Cliffordian from\ $\R\oplus V_n$\ to\
$\R_{0,n}$.
\bigskip
{\bf Proof.-}
First we show that if\ $x\in\R\oplus V_n$\ then\ $P_\beta (x)\in\R_{0,n}$\ or
better\ $P_\beta (x)\in\R\oplus V_n$. We know that\ $P_\beta (z) = D_\ast
p_\beta (z)$\ in\ $\R_{0,n+1}$\ and\ $p_\beta (z)$\ is even in\ $z_{n+1}$\
since\ $\beta_{n+1}$\ is even. So $\left[\displaystyle{\partial\over \partial
z_{n+1}} p_\beta (z)\right]_{z_{n+1}=0} = 0$\ and\ $P_\beta (x)\in \R\oplus
V_n$\ for\ $x = {\cal R}z$. 

Secondly we show that\ $P_\beta (x)$\ is analytic Cliffordian. We know that 
$$P_\beta (x) = \sum_{\sigma\in S_\beta} \Bigl(\prod_{\nu = 1}^n
(e_{\sigma(\nu )}x)\Bigr)\ e_{\sigma (n+1)},$$

which means that\ $P_\beta (x)$\ is the sum of all different polynomials
deduced from

$$\underbrace{ e_0xe_0\ldots e_0}_{\beta_0\hbox{times}\, e_0} x 
\underbrace{ e_1xe_1\ldots e_1}_{\beta_1\hbox{times}\, e_1}
xe_2xe_2\ldots\ldots e_n x
\underbrace{ e_{n+1}x e_{n+1} x\ldots x_n e_{n+1}}_{\beta_{n+1}\hbox{times}\,
e_{n+1}}\leqno (\ast)$$

by permutations of the\ $e_i's$. Note that

$$e_{n+1} = (-1)^m e_{1\, 2\ldots n} e_{1\, 2\ldots n+1}$$

and that the pseudo scalar\ $e_{1\,2\ldots 2m\, 2m+1}$\ belongs to the center
of\ $\R_{0,2m+1}$. Since\ $\beta_{n+1}$\ is even and since \ $(e_{1\, 2\ldots
n+1})^2\in \{ 1,-1\}$, we deduce that up to a sign we can replace ($\ast$) by
:

$$\underbrace{ e_0xe_0\ldots e_0}_{\beta_0\hbox{times}\, e_0} x 
\underbrace{ e_1xe_1\ldots e_1}_{\beta_1\hbox{times}\, e_1}
xe_2xe_2\ldots\ldots e_n x
\underbrace{ e_{1\, 2\ldots n}x e_{1\, 2\ldots n} \ldots x 
e_{1\, 2\ldots n}}_{\beta_{n+1}\hbox{times}\, e_{1\, 2\ldots n}}\leqno
(\ast\ast)$$

Let us chose a basis in\ $V_n$\ such that\ $x = x_0e_0+x_1e_1$. Then\ $e_1$\
commutes with\ $x$\ and for\ $i\geq 2$, $e_ie_{i+1}$\ commutes with\ $e_1$,
$x$ and of course\ $e_0$. Suppose\ $\beta_2 > 0$, for each term of the form
$$Ae_2 Be_{1\, 2\ldots n} C$$
we have another term equal to\ $Ae_{1\, 2\ldots n}Be_n C$. But then the
commutation  rules give since\ $n = 2m$ 
$$Ae_2Be_{1\, 2\ldots n}C + Ae_{1\, e\ldots n}Be_2 C = 0.$$

Thus\ $\beta_2 = 0$. And similarly\ $\beta_3=\ldots =\beta_n = 0$. We get :
$$P_\beta (x) = \sum_{\sigma \in S_\beta} \Bigl(\prod_{\nu =1}^{|\beta |-1}
E_{\sigma (\nu )}x\Bigr) E_{\sigma (|\beta |)}$$

where\ $E_0 = e_0$, $E_1 = e_1$\ and\ $E_{n+1} = e_{1\, 2\ldots n}$. Commuting
systematically\ $e_1$\ and\ $e_3e_4, e_5e_6,\ldots ,e_{n-1}e_n$\ from\
$E_{n+1}$ to the right 
$$P_\beta (x) = \Bigl( \sum_{\sigma \in S_\alpha} \Bigl(\prod_{\nu =
1}^{|\alpha |-1} e_{\sigma (\nu )}x\Bigr)e_{\sigma (|\alpha |)}\Bigr)
(e_1)^{\beta_{n+1}}(e_{3\, 4\ldots n})^{\beta_{n+1}}$$

where\ $\alpha = (\beta_0,\beta_1,\beta_{n+1},0,\ldots ,0)\in (\{
0\}\cup\N)^{n+1}$. Thus\ $P_\beta (x) = \pm P_\alpha (x)$\ with\ $P_\alpha$\
analytic Cliffordian.

\bigskip
{\bf Lemma 2.-} Let\ $n = 2m$\ and\ $\beta = (\beta_0,\ldots ,\beta_n,
\beta_{n+1})$. If\ $\beta_{n+1}$\ is odd then the restriction of\ $P_\beta
e_{n+1}$\ to\ $\R\oplus V_n$\ is analytic Cliffordian from\ $\R\oplus V_n$\ to\
$\R_{0,n}$.
\medskip
{\bf Proof.-}
Now\ $p_\beta (z)$\ is odd in\ $z_{n+1}$\ so\
$\left[\displaystyle{\partial\over \partial z_i} p_\beta
(z)\right]_{z_{n+1}=0} = 0$\ for\ $i = 0,\ldots , n$\ and\ $P_\beta (x) = -
\left[\displaystyle{\partial\over \partial z_{n+1}} p_\beta
(z)\right]_{z_{n+1} =0} e_{n+1}$\ and\ $P_\beta (x) e_{n+1}$\ is a scalar. For
the second part of the proof we reduce the sum

$$P_\beta (x) e_{n+1} = \sum_{\sigma\in S_\beta} \Bigl(\prod_{\nu=1}^{|\beta
|-1} e_{\sigma (\nu )}x\Bigr) e_{\sigma (|\beta |)} e_{n+1}$$
 as in Lemma 1.

\bigskip\bigskip
{\bf Corollary 1.-} \ The space of left analytic Cliffordian functions is an \
$\R$-vector space and an \ $\R_{0,n}$-right module, closed relatively to scalar
derivations.

\bigskip\medskip
To deduce corollary 2 and 3 from the above theorem, the following lemma is convenient.

\smallskip
{\bf Lemma 3.-} \ If \ $v\in\R\oplus V_{2m}$, \ then \ $\displaystyle\sum_{i=0}^{2m}
e_i ve_i = (1-2m) v_*$.

\bigskip\bigskip
{\bf Corollary 2.-} \ If \ $f$ \ is left analytic Cliffordian, then \ $Df$ \ is also
left analytic Cliffordian.
\medskip
{\bf Proof.-} If \ $n$ \ is odd let \ $n = 2m+1$. Since \ $f$ \ is left analytic
Cliffordian it is left holomorphic Cliffordian and \ $D\Delta^mf = 0$. But since \ $D$
\ commutes with \ $D$ \ and \ $\Delta$ \ we have
$$D\Delta^m (Df) = D(D\Delta^mf) = 0.$$

Thus \ $Df$ \ is left holomorphic Cliffordian or left analytic Cliffordian.

If \ $n$ \ is even let \ $n = 2m$.  Let us write ~: $x = z_0+z_1e_1 +\ldots+z_{2m}
e_{2m}$ \ and
\ $D' = \displaystyle\sum_{i=0}^{2m} \ e_i \ \displaystyle{\partial\over \partial z_i}.
\quad
\hbox{Then we have}\quad D = D' + e_{2m+1} \ \displaystyle{\partial\over \partial
z_{2m+1}}.$ 

We want to show that if \ $u$ \ is left analytic Cliffordian then  \ $D'u$
\ is also left analytic Cliffordian. In fact, we need only to show this for \ $u(x) =
M_N^a(x)$
\ for any \ $a\in\R \oplus V_{2m}$ \ and any \ $N\in\N$.
A straightforward computation gives
$$D'M_N^a (x) = \sum_{k=1}^{N-1} \ \sum_{i=0}^{2m} \ e_i \ M_k^a (x) e_i \ M_{N-k}^a
(x)$$ and since \ $M_k^a (x) \in \R\oplus V_{2m}$, lemma 3 gives us
$$D'M_N^a (x) = - (2m-1) \sum_{k=1}^{N-1} {\left[ M_k^a (x)\right]}_* \ M_{N-k}^a
(x).$$
If \ $N$ \ is odd, let \ $N = 2M+1$; \ we have
$$D'M_{2M+1}^a (x) = -2(2m-1) \sum_{k=1}^M \ \vmid{a}^{2k} \ \vmid{x}^{2k-2}
S\bigl( (xa)^{2M-2k+1}\bigr).$$

\bigskip
If \ $N$ \ is even, let \ $N = 2M$; \ we have
$$D'M_{2M}^a (x) = - (2m-1) \ \left\{ \vmid{a}^{2M} \ \vmid{x}^{2M-2} + 2
\sum_{k=1}^M  \ \vmid{a}^{2k} \ \vmid{x}^{2k-2} \  S\bigl(
(xa)^{2M-2k}\bigr)\right\}.$$

\bigskip
The same computations in \ $\R_{0,2m+1}$ \ give
$$D~\!M_{2M+1}^a(z) = -2(2m) \ \sum_{k=1}^M \ \vmid{a}^{2k} \ \vmid{z}^{2k-2} \
S\bigl( (za)^{2M-2k+1}\bigr)$$
and
$$D~\!M_{2M}^a(z) = - (2m) \ \left\{ \vmid{a}^{2M} \ \vmid{z}^{2M-2} +
2~\sum_{k=1}^M \ \vmid{a}^{2k} \ \vmid{z}^{2k-2} \ S\bigl(
(za)^{2M-2k}\bigr)\right\}.$$

\bigskip
Then we have for any \ $N$
$$D' M_N^a (x) = {2m-1\over 2m} \ {\left[D~\! M_N^a (z)\right]}_{z=x}.$$

Since \ $a\in\R_{0,2m}$, \ we have \ $S\bigl( (za)^{N-2k}\bigr) = S\bigl(
(xa)^{N-2k}\bigr)$, and since \ $\vmid{\overline z}~=~\vmid{z}$, \ we have \ $D~\!
M_N^a (\overline z) = D~\!\! M_N^a (z) = \overline{D~\!\! M_N^a(z)}$.  The theorem
then proves that \ $D'M_N^a(x)$ \ is left analytic Cliffordian.

\bigskip\bigskip
{\bf Corollary 3.-} \ If \ $f$ \ is left analytic Cliffordian, then \ $D_* f_*$ \ is
also left analytic Cliffordian.
\medskip
{\bf Proof.-} Let us write \ $f$ \ as 
$$f(z) = \sum_{N=1}^\infty \ \sum_{\vmid{\ \!\alpha\ \!}= N} P_\alpha (z) c_\alpha =
\sum_{N=1}^\infty \ \sum_{\vmid{\ \!\alpha\ \! }= N} D_*~\!p_\alpha (z) c_\alpha$$
and define the left analytic Cliffordian function \ $\widetilde f$ \ by
$$\widetilde f(z) = \sum_{N=1}^\infty ~\sum_{\vmid{~\!\alpha~\!}= N} \ P_\alpha(z)
c_{\alpha *}.$$

Then we have
$$D_* \ f_*(z) = {\bigl( Df(z)\bigr)}_* = \sum_{N=1}^\infty \
\sum_{\vmid{~\!\alpha~\!}=N} \ \Delta p_\alpha (z) \ c_{\alpha *} = D~\! \widetilde
f(z).$$

By corollary 1, \ $D\widetilde f$ \ is left analytic Cliffordian, thus \ $D_* f_*$ \
is also left analytic  Cliffordian.

\vskip1cm
\hskip-0,9cm{\bf  V. Cauchy's problem and boundary data}
\bigskip
{\bf  V.1. Aim} 
\medskip
We intend to generalize to Clifford algebras of any dimensions the method of
extending real analytic functions into complex holomorphic functions. Let  \
$\Omega$ \ be a domain of \ $\R \oplus V_{2m}$ \ and \ $u : \Omega \fle\R_{0,2m}$.
The Cauchy problem
$$\cases{D\Delta^mf = 0\cr
\quad\ \  f\mid_\Omega = u\cr}$$

where the unknown function \ $f$ \ has to be defined on an open set \ $\Lambda$ \
containing \ $\Omega$, \ seems not well defined since the partial differential
equations are of order $2m+1$. The Cauchy-Kowalewski theorem tells us that we need
the normal derivatives of \ $f$ \  to \ $\R\oplus V_{2m}$ \ up to the order $2m$. We
will see that the algebraic and analytic properties of \ $f$ \ and \ $u$ \ enable us
to compute these derivatives uniquely given the function $u$.
\medskip
We will denote by \ ${\cal A}_{2m}$ \ the linear space of left analytic Cliffordian
functions defined on \ $\Omega$ \ and taking their values in \ $\R_{0,2m}$.
Similarly, \ ${\cal A}_{2m+1}$ \ is the linear space of left analytic Cliffordian
functions defined on \ $\Lambda$ \ and with values in \ $\R_{0,2m+1}$.

\bigskip\bigskip
{\bf  V.2. The operator \ $(A\mid \nabla_n)$} 
\bigskip
{\bf Definition.-}  Let \ $A\in\R_{0,n}$.  We define \ $(A\mid\nabla_n) :
{\cal A}_n \flep {\cal A}_n$ \ by:
\smallskip
{\parindent=1,5cm
\item{(i)} $\forall a\in\R\oplus V_n \hskip0,3cm \forall N\in\N \qquad
(A\mid\nabla_n) M_N^a (z) = \displaystyle\sum_{k=1}^{N-1} \ M_k^a (z) \ A \
M_{N-k}^a (z)$
\medskip
\item{(ii)} $\forall f\in {\cal A}_n \hskip0,9cm \forall K\in\R_{0,n} \ \ 
(A\mid\nabla_n) (f(z)K) = \bigl( A\mid\nabla_n) f(z)\bigr) K$
\bigskip
\item{(iii)} $(A\mid\nabla_n)$ \ is \ $\R$-linear and continuous.
\par}

\bigskip
{\bf Consequence.-}  If \ $f(z) = \displaystyle\sum_{N=1}^\infty \ \sum_{a\in A_n} \
M_N^a (z) C_a$, \ then
$$(A\mid\nabla_n) f(z) = \sum_{N=1}^\infty \ \sum_{a\in A_N} \ \bigl(
(A\mid\nabla_n) M_N^a (z)\bigr) C_a.$$

\bigskip
{\bf Proposition 9.-} \ If \ $f\in {\cal A}_n$, \ then \
$\displaystyle{\partial\over\partial z_k} f(z) = (e_k\mid\nabla_n) f(z)$.
\medskip
{\bf Proof.-} We need only to verify the proposition for \ $f(z) = M_N^a(z)$. We
have
\vskip-0,6cm
$$\hskip-0,7cm{\partial\over\partial z_k} azaz \ldots za = a~e_k az \ldots za + azae_k
\ldots za +\ldots+ azaz \ldots e_k a = (e_k \mid \nabla_n) azaz \ldots za.$$
\smallskip

\bigskip\bigskip
{\bf Proposition 10.-} If \ $f\in {\cal A}_n$ \ and if \ $\lambda$ \ belongs to the
center of \ $\R_{0,n}$ \ then
$$(\lambda A \mid \nabla_n) f = \lambda (A\mid\nabla_n)f.$$

\medskip
{\bf Lemma.-} If \ $v\in\R\oplus V_{2m}$, \ then
$$\sum_{i=0}^{2m} \ e_i v e_i = (-1)^m (1-2m) e_{1\,2\ldots 2m} \ v e_{1\,2\ldots
2m}.$$
\medskip
{\bf Proof.-} Both sides of the equality are equal to \ $(1-2m)v_*$.

\bigskip\bigskip
{\bf Proposition 11.-} If \ $u\in {\cal A}_{2m}$, \ then
$$e_{1\,2\ldots 2m} (e_{1\,2\ldots 2m} \mid \nabla_{2m}) u = {(-1)^{m+1}\over
2m-1} Du.$$
\medskip
{\bf Proof.-} We need only verify the relation for \ $u(x) = M_N^a(x)$. Using the
definition and the lemma, we get
$$\eqalign{
(-1)^m (1-2m) e_{1\,2\ldots 2m} (e_{1\,2\ldots 2m} \mid \nabla_{2m}) u &=
\sum_{k=1}^{N-1} \ \sum_{i=0}^{2m} \ e_i M_k^a (x) e_i \ M_{N-k}^a (x)\cr
&= \sum_{i=0}^{2m} \ e_i \ \sum_{k=1}^{N-1} \ M_k^a (x) e_i M_{N-k}^a (x)\cr
&= \sum_{i=0}^{2m} \ e_i \ {\partial\over\partial x_i} \ M_N^a(x) = D \ \! M_N^a
(x).\cr}$$

\bigskip
{\bf Proposition 12.-} For \ $\mu\in \{ 0\} \cup \N$, \ we have
$$\eqalign{
&(e_{1\,2\ldots 2m} \mid \nabla_{2m})^{2\mu} u = {(-1)^{m\mu}\over
(2m-1)^{2\mu}} \ \Delta^\mu u\cr
&(e_{1\,2\ldots 2m} \mid \nabla_{2m})^{2\mu +1}u = {(-1)^{1+m\mu}\over
(2m-1)^{2\mu+1}} e_{1\,2\ldots 2m} \ D\Delta^\mu u.\cr}$$
\medskip
{\bf Proof.-} From proposition 3, we get
$$\eqalign{
(e_{1\,2\ldots 2m} \mid \nabla_{2m}) u &= {-1\over 2m-1} \ e_{1\,2\ldots 2m}\ 
Du\cr &= {-1\over 2m-1} \ D_* \ u_* e_{1\,2\ldots 2m}\cr}$$

The corollary 2 of the theorem of last paragraph shows that \ $D_*u_*e_{1\,2\ldots
2m}$ \ is left analytic Cliffordian. Thus we may apply the operator as many times as
we want. We get
$$\eqalign{\bigl( (e_{1\,2\ldots 2m} \mid \nabla_{2m}\bigr)^2 u
&= {-1\over 2m-1} \ e_{1\,2\ldots 2m} \ D \left\{ {-1\over 2m-1} D_* u_* e_{1\,2\ldots
2m} \right\}\cr 
&= {1\over (2m-1)^2} \ e_{1\,2\ldots 2m} \ \Delta u_* \ e_{1\,2\ldots 2m} \cr
&= {1\over (2m-1)^2} (e_{1\,2\ldots 2m})^2 \ \Delta u\cr
&= {(-1)^m\over (2m-1)^2} \ \Delta u. \cr}$$

A simple recursion gives then the general formulae.

\vfill\eject
{\bf  V.3. Analytic extention}
\bigskip
{\bf Theorem.-} \ Let \ $\Omega$ \ be a domain of \ $\R \oplus V_{2m}$ \ and let \
$u : \Omega \fle \R_{0,2m}$ \ be left analytic Cliffordian. There exist a domain \
$\Lambda$ \ in \ $\R\oplus V_{2m+1}$ \ with \ $\Omega\subset\Lambda$
\ and a unique left holomorphic Cliffordian fonction \ $f$ \ defined on \ $\Lambda
\subset
\R\oplus V_{2m+1}$,  \ such that \ $f\mid_\Omega =
u$. That function \ $f$ \  is such that \ $f(\overline z) =
\overline{f(z)}$ \ and if we denote by \ $\displaystyle{\partial\over\partial n}$ \
the normal derivative  to \ $\Omega$, \ we have for any \ $\mu$ \ in \ $\{ 0\} \cup
\N$
$$\eqalign{
&{\bigl( {\partial\over \partial n}\bigr)}^{2\mu} f\mid_\Omega = {(-1)^\mu\over
(2m-1)^{2\mu}} \Delta^\mu u\cr 
&{\bigl( {\partial\over \partial n}\bigr)}^{2\mu +1} \ f\mid_\Omega =
{(-1)^{\mu +1}\over (2m-1)^{2\mu+1}} \ e_{2m+1} \ D\Delta^\mu~u.\cr}$$
\medskip
{\bf Proof.-} {\bf 1$^\circ$)} Suppose \ $f$ \ is a solution. We use the usual notations
$$z = z_0+z_1 e_1 +\ldots +z_{2m} e_{2m} + z_{2m+1} e_{2m+1} = x + z_{2m+1}
e_{2m+1}.$$
Thus we have 
$${\bigl({\partial\over \partial n}\bigr)}^j f = \bigl( {\partial\over \partial
z_{2m+1}}\bigr)^j \ f.$$
Since \ $f$ \ is left analytic Cliffordian, we have \ $\displaystyle{\partial\over
\partial z_{2m+1}} f = \bigl( e_{2m+1} \mid \nabla_{2m+1}\bigr) f$ \ and thus
$${\bigl( {\partial\over\partial n}\bigr)}^j \ f = \bigl( e_{2m+1} \mid
\nabla_{2m+1}\bigr)^j f.$$

Note that \ $e_{2m+1} = (-1)^m \ e_{1\, 2\ldots 2m \ 2m+1} \  e_{1\,2\ldots 2m}$  \
and that \ $(-1)^m \ e_{1\, 2\ldots 2m \ 2m+1}$ \ belongs to the center of \
$\R_{0,2m+1}$. Using proposition 10 we get
$${\bigl( {\partial\over \partial n}\bigr)}^j \ f = \bigl( (-1)^m \ e_{1 \, 2\ldots
2m \ 2m+1}\bigr)^j \ (e_{1\, 2\ldots 2m} \mid \nabla_{2m+1})^j ~f.
$$
\medskip
Taking the restriction to \ $z = x\in\Omega$, \ we get
$${\bigl( {\partial\over \partial n}\bigr)}^j \ f\mid_\Omega (x) = (-1)^{mj} \
(e_{1\, 2\ldots 2m \ 2m+1})^j \ (e_{1\, 2\ldots 2m} \mid \nabla_{2m})^j \
u(x).$$

\medskip
Proposition 12 gives us then for \ $j = 2\mu$ \ and \ $j = 2\mu + 1$
$${\bigl( {\partial\over \partial n}\bigr)}^{2\mu} \ f\mid_\Omega (x) = {(-1)^\mu
\over (2m-1)^{2\mu}} \ \Delta^\mu u(x)$$
and
$${\bigl( {\partial\over \partial n}\bigr)}^{2\mu +1} \ f\mid_\Omega (x) =
{(-1)^{\mu +1}\over (2m-1)^{2\mu +1}} \ e_{2m+1} \ D\Delta^\mu u(x).$$

\medskip
{\bf 2$^\circ$)} Using the above conditions for \ $j\leq 2m$, \ the theorem of
Cauchy-Kowalewski proves the existence of \ $\Lambda'\subset\R\oplus V_{2m+1}$ \ with
\ $\Lambda'\cap (\R\oplus V_{2m}) = \Omega$ \ and the existence and unicity of \ $f$
\ in \ $\Lambda'$.

\medskip
{\bf 3$^\circ$)} Knowing the existence of $f$, the usual Taylor formula gives us
$$\eqalign{
f(z) = &\sum_{\mu = 0}^\infty \ {(-1)^\mu \over (2m-1)^{2\mu} (2\mu)!} \
(z_{2m+1})^{2\mu}
\ \Delta^\mu \ u(x)\cr
&\qquad + e_{2m+1} \ \sum_{\mu=0}^\infty \ { (-1)^{\mu +1}\over (2m-1)^{2\mu +1}
(2\mu +1)!} \ (z_{2m+1})^{2\mu +1} \ D\Delta^\mu u(x).\cr}$$
\medskip
This formula shows, since \ $Du(x) \in\R_{0,2m}$, \ that \ $f(\overline z) =
\overline{f(z)}$ \ in a subdomain \ $\Lambda$ \ of \ $\R\oplus V_{2m+1}$ \ such that
\ $\Omega\subset\Lambda\subset\Lambda'$.

\vskip1cm
\hskip-0,9cm{\bf  VI. Fueter's method}
\bigskip
Fueter's method is well known and
widely used to construct functions connected to monogenic functions \ [Fu], [De],
[Qi]. It is known to be effective to construct holomorphic Cliffordian functions in
the case of odd $n$. We show that it is still valid for our definition in the case of
even \ $n$.

\bigskip
{\bf Theorem.-}  Let \ $\varphi$ \ be a complex holomorphic function defined on \
$D_\varphi$ \ an open subset of the upper half-plane and let \ $p$ \ and \ $q$ \ be the
real functions of two variables defined by
$$\forall\zeta = \xi + i\eta \in D_\varphi \qquad \varphi (\zeta ) = p (\xi,\eta)
+ i q(\xi, \eta).$$
Let \ $\overrightarrow z = z_1 e_1 +\ldots + z_ne_n$, \ for \ $z_0+i \
\vmid{\overrightarrow z}\in D_\varphi$ \ we define \ $u(z_0+\overrightarrow z)$ \ by~:
$$u (z_0 + \overrightarrow z) = p (z_0, \ \vmid{\overrightarrow z}) +
{\overrightarrow z\over \vmid{\overrightarrow z}} q (z_0, \ \vmid{\overrightarrow
z}).$$

Then \ $u$ \ is a (left and right) holomorphic Cliffordian function.
\medskip
{\bf Proof.-} For odd \ $n$, this result is already known [LR1].

If \ $n$ \ is even, let \ $n = 2m$. Let \ $x = z_0 + \overrightarrow z$ \ and \ $z =
x + z_{2m+1} \ e_{2m+1}$. Define \ $f$ \ by \ $f(z) = p(z_0, \ \vmid{\overrightarrow
z + z_{2m+1} \ \! e_{2m+1}}) + \displaystyle{\overrightarrow z + z_{2m+1} \
e_{2m+1}\over
\vmid{\overrightarrow z + z_{2m+1} \ e_{2m+1}}}~q(z_0, \ \vmid{\overrightarrow z +
z_{2m+1} \ e_{2m+1}})$.  
\smallskip 
From the case of odd \ $n$, we know that \ $f$ \ is a
left and right holomorphic Cliffordian function. The theorem of paragraphe IV shows
then that
\ $u$ \ is a left and right holomorphic Cliffordian
 function, since we have \
$f(\overline z) =
\overline{f(z)}$ \ and \ $u(x) = f(x)$.

\bigskip\bigskip
\centerline{\bf Bibliographie}
\bigskip\bigskip
\bb BDS&{\sc F. Bracks, R. Delanghe, F. Sommen}&Clifford analysis&Pitman (1982) & &

\bb De&{\sc C. Deavors}&The quaternion calculus&Am. Math. monthly (1973), 995-1008
& &

\bb DSS&{\sc  R. Delanghe, F. Sommen, V. Soucek}&Clifford algebra and
spinor-valued functions&Kluwer (1992)& &

\bb EL1&{\sc S.L. Eriksson-Bique, H. Leutwiler}&On modified quaternionic analysis
in $\R^3$&Arch. Math.70 (1998), 228-234& &

\bb EL2&{\sc S.L. Eriksson-Bique, H. Leutwiler}&Hyperholomorphic functions&to appear&
&

\bb Fu&{\sc R. Fueter}&Die Functionentheorie der Differengleichungen  $\Delta u
= 0$  and  $\Delta\Delta u = 0$ mit vier reellen Variablen&Comm. Math. Helv.7 (1935),
307-330& &

\bb Le1&{\sc H. Leutwiler}&Modified quaternionic analysis in $\R^3$&Complex
variables Vol.20 (1992), 19-51& &

\bb Le2&{\sc H. Leutwiler}&Rudiments of a function theory in $\R^3$&Expo. Math.$\!\!$
14
 (1996), 97-123& &

\bb LR1&{\sc G. Laville, I. Ramadanoff}&Holomorphic Cliffordian functions&Advances
in Applied Clifford algebras 8, n$^\circ$2 (1998), 321-340& &

\bb LR2&{\sc G. Laville, I. Ramadanoff}&Elliptic Cliffordian functions&Complex
variables Vol. 45, n$^\circ$4 (2001), 297-318& &

\bb Qi&{\sc T. Qian}&Generalization of Fueter's result to $\R^{n+1}$&Rend. Math.
Acc. Lincei 8 (1997), 111-117& &

\end